\documentclass[11pt,a4paper]{article}
\usepackage[utf8]{inputenc}
\usepackage[english]{babel}%remembertwochangestoactivate
\usepackage{amsfonts} 
\usepackage{amsmath}
\usepackage{amssymb}
\usepackage{amsthm}%%%Proof environment
\usepackage{mathrsfs}%%%Mathscr riccioli
\usepackage{microtype}%%%secretpackagemarcello
\usepackage{marvosym}%%%always BEFORE wasysym or arxiv crashes
\usepackage{wasysym} %%%diameter=varnothing=emptyset
\usepackage{stmaryrd}
\usepackage{dsfont}%%% alternative to previous 1, 0
\usepackage{enumerate}
\usepackage{multicol,multirow}%%%column(s)
\usepackage{graphicx}%%%pics
\usepackage{xcolor}
\usepackage{authblk}%%%multiauthors
\usepackage{url}%%%size of urls. default was {\scriptsize}

\usepackage{hyperref}%%%linklink
\hypersetup{
	colorlinks=false,
	linkcolor=blue,
	urlcolor=red,
	linktoc=all}

\usepackage{tikz}%%%TikZ
\usetikzlibrary{decorations.pathreplacing}
\usetikzlibrary{decorations.pathmorphing}
\usetikzlibrary{calc}
\usetikzlibrary{decorations.pathreplacing}
\usetikzlibrary{decorations.pathmorphing}
\usetikzlibrary{arrows}
\usetikzlibrary{patterns,fadings}
\usepackage[all,cmtip]{xy}%%%xy alternative to tikz-cd

\usepackage[tocflat]{tocstyle}%alpha version, dangerous
\usepackage{tocloft}
\newcommand{\Hilb}{\mathcal{H}}%myHilb
\newcommand{\oneop}{\mathds{1}}
\newcommand{\iotabar}{{\overline{\iota}}}

\newcommand{\Z}{\mathcal{Z}}

\newcommand{\A}{\mathcal{A}}
\newcommand{\B}{\mathcal{B}}

\newcommand{\cL}{\mathcal{L}}
\renewcommand{\L}{\cL}%rinominato comando \L
\newcommand{\M}{\mathcal{M}}

\newcommand{\N}{\mathcal{N}}

\newcommand{\R}{\mathcal{R}}

\renewcommand{\L}{\cL}

\newcommand{\RR}{\mathbb{R}}

\newcommand{\CC}{\mathbb{C}}

\newcommand{\NN}{\mathbb{N}}

\DeclareMathOperator{\End}{End}
\DeclareMathOperator{\Bim}{Bim}
\DeclareMathOperator{\Ind}{Ind}

\DeclareMathOperator{\id}{id}

\DeclareMathOperator{\tr}{tr}

\def\III{{I\!I\!I}}%typeIII khr
\def\II{{I\!I}}%typeII khr

\newcommand{\lqq}{\lq\lq}

\usepackage{xspace}
\DeclareRobustCommand{\eg}{e.g.\@\xspace}
\DeclareRobustCommand{\cf}{cf.\@\xspace}
\DeclareRobustCommand{\ie}{i.e.\@\xspace}

\DeclareRobustCommand{\Sec}{Sec.\@\xspace}
\DeclareRobustCommand{\Prop}{Prop.\@\xspace}
\DeclareRobustCommand{\Lem}{Lem.\@\xspace}
\DeclareRobustCommand{\Cor}{Cor.\@\xspace}
\DeclareRobustCommand{\Thm}{Thm.\@\xspace}
\DeclareRobustCommand{\Ch}{Ch.\@\xspace}

\DeclareRobustCommand{\wrt}{w.r.t.\@\xspace}
\makeatletter
\DeclareRobustCommand{\etc}{%
    \@ifnextchar{.}%
        {etc}%
        {etc.\@\xspace}%
}
\makeatother

\newcommand{\Cstar}{$C^\ast$\@\xspace}

\def\u1net{{\A_\RR}}

\theoremstyle{plain}
\newtheorem{theorem}{Theorem}[section]

\newtheorem{proposition}[theorem]{Proposition}

\newtheorem{problem}[theorem]{Problem}

\theoremstyle{definition}
\newtheorem{definition}[theorem]{Definition}

\theoremstyle{remark}

\newtheorem{remark}[theorem]{Remark}

\numberwithin{equation}{section}%number equations with section in amsmath

\begin{document}

\title{\LARGE Minimal index and dimension for inclusions of von Neumann algebras with finite-dimensional centers}

\author{\Large Luca Giorgetti\thanks{\footnotesize Supported by the ERC Advanced Grant n.\ 669240 QUEST ``Quantum Algebraic Structures and Models'', 
%MIUR FARE R16X5RB55W QUEST-NET, 
OPAL ``Consolidate the Foundations" and GNAMPA--INdAM.}
}
%\author{\Large Roberto Longo}
%%%\author[1]{Author \thanks{E-mail}}
%\affil{\normalsize Dipartimento di Matematica, Universit\`a di Roma Tor Vergata\\

%Via della Ricerca Scientifica, 1, I-00133 Roma, Italy\\

%{\tt giorgett@mat.uniroma2.it}, {\tt longo@mat.uniroma2.it}}

\affil{\normalsize Dipartimento di Matematica \lqq Guido Castelnuovo"\\ 

Sapienza Universit\`a di Roma\\ 

Piazzale Aldo Moro, 5, I-00185 Roma, Italy\\

{\tt giorgetti@mat.uniroma1.it}}

%\affil{\normalsize Dipartimento di Matematica, Universit\`a di Roma Tor Vergata\\

%Via della Ricerca Scientifica, 1, I-00133 Roma, Italy\\

%{\tt giorgett@mat.uniroma2.it}, {\tt longo@mat.uniroma2.it}}
%%%nolastcomma\renewcommand\Authands{ and }

\date{}

\maketitle

\begin{abstract}
The notion of index for inclusions of von Neumann algebras goes back to a seminal work of Jones on subfactors of type $\II_1$. In the absence of a trace, one can still define the index of a conditional expectation associated to a subfactor and look for expectations that minimize the index. This value is called the minimal index of the subfactor. 

We report on our analysis, contained in \cite{GiLo19}, of the minimal index for inclusions of arbitrary von Neumann algebras (not necessarily finite, nor factorial) with finite-dimensional centers. Our results generalize some aspects of the Jones index for multi-matrix inclusions (finite direct sums of matrix algebras), e.g., the minimal index always equals the squared norm of a matrix, that we call \emph{matrix dimension}, as it is the case for multi-matrices with respect to the Bratteli inclusion matrix. We also mention how the theory of minimal index can be formulated in the purely algebraic context of rigid 2-{$C^{\ast}$}-categories. 
%Joint work with Roberto Longo (Roma Tor Vergata), preprint arxiv:1805.09234.
\end{abstract}

%\tableofcontents

%\vskip3cm

%{\footnotesize Supported by the ERC Advanced Grant 669240 QUEST ``Quantum Algebraic Structures and Models'', 
%MIUR FARE R16X5RB55W QUEST-NET, 
%OPAL\,-``Consolidate the Foundations", GNAMPA-INdAM.}

\section{Motivation}

One motivation for studying von Neumann algebras with non-trivial centers and inclusions, or better bimodules, between them comes from the theory of \emph{Quantum Information}. 

In an operator-algebraic description of quantum systems, \emph{observables} are described by the self-adjoint part of a non-commutative von Neumann algebra $\M$ (with separable predual), while \emph{states} correspond to normal faithful positive functionals $\varphi:\M\rightarrow\CC$ normalized such that $\varphi(\oneop) = 1$, where $\oneop$ denotes the identity operator. Keep in mind as an example the most commonly studied case of \emph{finite} quantum systems \cite[Part I]{OhPeBook} where the algebra generated by the observables is finite-dimensional, thus a multi-matrix algebra. Namely, $\M\cong \bigoplus_{i=1,\ldots,m} M_{k_i}(\CC)$, where $m, k_i \in \NN$ and $M_{k_i}(\CC)$ is the algebra of $k_i\times k_i$ matrices over $\CC$, realized on the finite dimensional Hilbert space $\CC^N$, $N = k_1 + \ldots + k_m$. More generally, the center $\Z(\M) = \M \cap \M'$ is the classical part of the system, in the previous case $\Z(\M) \cong \CC^m$, while each factor in the central decomposition of $\M$, in the previous case $M_{k_i}(\CC)$, is a purely quantum part of the system. Recall that a factor is a von Neumann algebra with center equal to $\CC\oneop$.

In this note we shall always assume, as in \cite{GiLo19}, that the center is \emph{finite-dimensional}, $\Z(\M) \cong \CC^m$, thus
$$\M \cong \bigoplus_{i=1,\ldots,m} \M_i$$
where $p_i\in\Z(\M)$ are the minimal central projections and $\M_i = \M p_i$ are factors (of arbitrary type). In this sense, we study possibly \emph{infinite} quantum systems with a finite-dimensional classical part.  

The building blocks of information transfer (communication) from a quantum system $\N$ to another $\M$ are called \emph{quantum channels}. In the operator-algebraic setting, they are conventionally described by normal completely positive maps $\alpha: \N \rightarrow \M$ such that $\alpha(\oneop) = \oneop$. Recall that a map $\alpha:\N\rightarrow\M $ is called positive if it sends positive elements of $\N$ to positive elements of $\M$, while it is called completely positive if $\alpha \otimes \id_{k\times k} : \N \otimes M_k(\CC) \rightarrow \M \otimes M_k(\CC)$ is positive for every $k\in\NN$. Communication takes place via transferring states from one system to another by pullback, namely $\alpha^t(\varphi) := \varphi \circ \alpha$ is a state on $\N$ whenever $\varphi$ is a state on $\M$. Note also that normal states and normal unital $^*$-homomorphisms are examples of completely positive maps. In this note, as in the first part of \cite{GiLo19}, we will mostly restrict ourselves to quantum channels given by inclusion morphisms 
$$\iota : \N \hookrightarrow \M$$
associated to unital inclusions of von Neumann algebras $\N\subseteq\M$. We will furthermore assume that the inclusion morphism has \emph{finite Jones index}. This assumption is equivalent to the existence of a \emph{conjugate} morphism $\iotabar : \M \rightarrow \N$ (conjugate quantum channel going in the opposite direction). The notion of Jones index will be reviewed in the next section, while we refer to \cite{LoRo97} and to the second part of \cite{GiLo19} for the definition of conjugate morphism and its relation to the theory of (minimal) index. For now, we only mention that these notions of conjugation and of minimal index, and the more fundamental notion of (matricial, intrinsic) dimension, are naturally formulated in a tensor \Cstar-categorical language. Namely for abstract 1-arrows $X:\N\rightarrow\M$ running between 0-objects $\N,\M$ of a \emph{2-\Cstar-category}.

\begin{remark}
More generally, one can think of quantum channels as described by $\N$-$\M$ bimodules $\Hilb$, also denoted by $_\N \Hilb _\M$, see \cite[\Sec 2,3]{Lon18}. Recall that a bimodule is a Hilbert space $\Hilb$ with a normal left action of $\N$, $l : \N \rightarrow \B(\Hilb)$, namely $l(n_1n_2) = l(n_1)l(n_2)$, and a normal right action of $\M$, $r : \M \rightarrow \B(\Hilb)$, namely $r(m_1m_2) = r(m_2)r(m_1)$, such that $l(\N)$ and $r(\M)$ mutually commute in $\B(\Hilb)$. Thus a bimodule sees in a \lqq balanced way" the inclusion $l(\N) \subset r(\M)'$ and the dual inclusion $r(\M) \subset l(\N)'$. Moreover, it is known that every normal unital completely positive map gives rise to a bimodule.
\end{remark}

In this algebraic setup of Quantum Information, Longo \cite[\Thm 3.2, \Cor 3.4]{Lon18} gave a mathematical derivation of \emph{Landauer's bound} for possibly infinite quantum systems \cite{Lan61}. See also \cite{Ben03}, \cite{PlVi01} for an introduction to Landauer's principle and bound, and for an explanation of how these settle the famous Maxwell's demon paradox. The bound is a lower estimate on the amount of energy (heat) that is emitted from the system whenever 1 classical bit of information is deleted (or any \emph{logically} irreversible operation is performed). Namely,
$$E_\alpha \geq \frac{1}{2} k T \log(2)$$
where $E_\alpha$ is the variation of the free energy associated to the channel $\alpha$, $k$ is Boltzmann's constant and $T$ is the temperature of the environment. The bound is calculated in \cite{Lon18} by means of the matrix dimension $D_\alpha$, it is in general half of the original Landauer's bound $E \geq k T \log(2)$, and it coincides with the latter in the case of \emph{finite} quantum systems because the lowest non-trivial possible (scalar) dimension of an inclusion of matrix algebras is $2$ instead of $\sqrt{2}$.

The most important properties of the matrix dimension $D_\alpha$ of a quantum channel $\alpha$, in our case of an inclusion morphism $\alpha = \iota$, are its \emph{multiplicativity} and \emph{additivity}:
$$D_{\beta \circ \alpha} = D_\beta D_\alpha, \quad D_{\alpha \oplus \beta} = D_\alpha + D_\beta$$
where $\beta \circ \alpha$ and $\alpha \oplus \beta$ denote respectively the composition (or \lqq tensor multiplication") and the direct sum of channels. Moreover, $D_\alpha$ determines the minimal index as the square of its $l^2$-norm, and the (unique) minimal conditional expectation via Perron-Frobenius theory, as we shall explain in the next section.

\section{Minimal index}

Let $\N \subset \M$ be a unital inclusion of von Neumann algebras. The \emph{Jones index} of the inclusion is a number ($\geq 1$) that measures the \lqq relative size" of $\M$ \wrt $\N$. The index equals $1$ if $\N = \M$, it equals $+\infty$ if $\M$ is way bigger than $\N$, and it enjoys the exciting property of being quantized between $1$ and $4$. Actually, there can be more than one notion of \lqq index" for an inclusion, depending on the type of algebras involved. If $\N \cong M_k(\CC)$ and $\M \cong M_h(\CC)$ (finite type $I$ factors) then $h=km$ for some $m\in\NN$, the inclusion morphism is the amplification, namely $M_k(\CC) \otimes \oneop_m \subset M_h(\CC)$, and the index is the square of the multiplicity $m^2$, \ie, the ratio of the algebraic dimensions of $\M$ over $\N$. If $\N$ and $\M$ are factors of type $\II_1$, the Jones index, denoted by $[\M:\N]$, is defined in terms of tracial states \cite{Jon83}. This is the original definition of index and in this regime one can observe the already mentioned quantization phenomenon of the index values \cite[\Thm 4.3.1]{Jon83}. If $\N \subset \M$ is a multi-matrix inclusion (the simplest instance of inclusion with non-trivial centers), namely $\N\cong \bigoplus_{j=1,\ldots,n} M_{k_j}(\CC)$ and $\M\cong \bigoplus_{i=1,\ldots,m} M_{h_i}(\CC)$, denote by $k = (k_1,\ldots,k_n)^t$, $h = (h_1,\ldots,h_m)^t$ the vectors of dimensions (with positive integer components), and by $\Lambda$ the Bratteli inclusion matrix (with positive integer entries) describing the inclusion morphism of $\N$ in $\M$, \cite{Bra72}. Then $\Lambda k = h$ is the only consistency condition on the Bratteli diagram associated to the inclusion, and the index (there are several equivalent definitions of index in this case) equals $\|\Lambda\|_{l^2}^2$, \cite[\Ch 2]{GdHJ89}. For an inclusion $\N\subset\M$ of finite direct sums of type $\II_1$ factors, the index, again denoted by $[\M:\N]$, is defined as the spectral radius of a product of matrices constructed from the unique trace on each factor in $\M$ (\emph{trace matrix}) and from the square roots of the Jones indices of the subfactors obtained by central decomposition (\emph{Jones index matrix}), \cite[\Ch 3]{GdHJ89}. If the inclusion is \emph{connected}, \ie, $\Z(\N)\cap\Z(\M) = \CC\oneop$ (which is equivalent to the actual connectedness of the adjacency graph of $\N\subset\M$, see \cite[\Sec 1.3]{GdHJ89}) and if $[\M:\N] \leq 4$, then $[\M:\N] = \|\Lambda^{\M}_\N\|_{l^2}^2$ by \cite[\Thm 3.7.13]{GdHJ89}, where $\Lambda^{\M}_\N$ is the aforementioned Jones index matrix. The theory of index for inclusion of finite von Neumann algebras has been further extended to inclusions with possibly infinite (atomic or diffuse) centers by Jolissaint in \cite{Jol90}.

In the absence of a trace, \eg, for subfactors of type $\III$, one can still define the index of a normal faithful \emph{conditional expectation} $E:\M\rightarrow\N$, denoted by $\Ind(E)$, \cite{Kos86}. An inclusion is said to have \emph{finite index} if it admits \emph{some} $E$ with $\|\Ind(E)\|<+\infty$. The index $\Ind(E)$ is in general an element of $\Z(\M)$, with $\Ind(E)\geq \oneop$, it is of course a scalar in the case of subfactors, and it gives back the Jones index for a finite subfactor by $[\M:\N]\oneop = \Ind(E^{\tr})$, where $E^{\tr}$ is the trace-preserving expectation. If the inclusion is not irreducible, \ie $\N'\cap\M \neq \CC\oneop$, there can be several expectations $E:\M\rightarrow\N$, and one can look for those \emph{minimizing} the number $\|\Ind(E)\|$. The \emph{minimal index} of the inclusion $\N\subset\M$ is then defined to be
$$[\M:\N]_0 := \inf_{E}\{\|\Ind(E)\|\}$$
This analysis has been performed first in the case of subfactors by Hiai \cite{Hia88}, Longo \cite{Lon89} and Havet \cite[\Sec 1]{Hav90}, where it is shown that there is a \emph{unique} expectation minimizing the index and this expectation is characterized via a certain \emph{sphericality condition} (which opens the way to a tensor \Cstar-categorical formulation of the minimal index, see \cite{LoRo97}).

\begin{remark}\label{rmk:extremality}
For finite subfactors it can happen that $[\M:\N]_0 \lneq [\M:\N]$. Equality is attained for the so called \emph{extremal} subfactors, see \cite[\Sec 4]{PiPo86}, \cite{PiPo91} and \cite[\Sec 2.3]{BurPhD}. 
%Note that inclusions of matrix algebras are always extremal, as well as irreducible or finite-depth type II_1 subfactors 
\end{remark}

In case of arbitrary inclusions, Jolissaint \cite[\Thm 1.8]{Jol91} proved that there is always an expectation minimizing the index (thus called a \emph{minimal expectation}), but this expectation is \emph{not} unique in general as shown by Fidaleo and Isola \cite[\Prop 10, \Sec 5]{FiIs96}.

The starting points of our analysis are the works of Havet \cite[\Sec 2]{Hav90} and Teruya \cite{Ter92}, where it is shown that in the case of \emph{connected} inclusions with \emph{finite-dimensional centers} there is a unique minimal expectation $E^0$, and its index is a scalar operator in $\Z(\M)$, \ie, $\Ind(E^0) = [\M:\N]_0 \oneop$. We mention that the connectedness assumption is almost without loss of generality, as every inclusion can be written as a direct sum of connected ones.

\begin{theorem}\emph{\cite{GiLo19}}.
Let $\N\subset\M$ be a connected inclusion with finite index and finite-dimensional centers. 

Let $p_1,\ldots,p_m$ and $q_1,\ldots,q_n$ be the minimal projections in $\Z(\M)$ and $\Z(\N)$ respectively and, whenever $p_iq_j \neq 0$, consider the subfactors $\N_{ij} := \N p_iq_j \subset \M_{ij} := q_j \M p_iq_j$ and define $d_{ij} := [\M_{ij}:\N_{ij}]^{1/2}_0$, while $d_{ij} := 0$ otherwise. Then 
$$[\M:\N]_0 = \| D\|_{l^2}^2$$
where $D$ is the $m\times n$ matrix with entries $d_{ij}$. We call $D$ the matrix dimension of $\N\subset\M$ and $d:=\| D\|_{l^2}$ its scalar dimension. 
\end{theorem}

From the previous theorem, together with the quantization of the index for expectations between factors \cite{Kos86}, as in \cite[\Prop 3.7.12 (c)]{GdHJ89} one can conclude that either $[\M:\N]_0\in\{4 \cos^2(\pi/k), k\in\NN, k\geq 3\}$ or $[\M:\N]_0 \geq 4$.

For every normal faithful conditional expectation $E:\M\rightarrow\N$ one can consider a matrix of expectations $E_{ij}:\M_{ij}\rightarrow\N_{ij}$ and a matrix of numbers $\lambda^E_{ij} \geq 0$ with the property that $\sum_{i}\lambda^E_{ij} = 1$ (thus called column-stochastic or Markovian). Namely, $\lambda^E_{ij} q_j := E(p_i)q_j$ and $E_{ij}(q_j x p_i q_j) := (\lambda^E_{ij})^{-1} E(x p_i) p_i q_j$ for every $x\in\M$.
From the pair $E_{ij}, \lambda^E_{ij}$ one can reconstruct the expectation via $E(x) = \sum_{i,j} \lambda^E_{ij} \sigma_{ij}(E_{ij}(q_j x p_i q_j))$, $x\in\M$, where $\sigma_{ij} : \N_{ij} \rightarrow \N_j := \N q_j$ is the inverse of the induction isomorphism $y q_j \mapsto y p_i q_j$, $y\in\N$. Moreover, every expectation arises uniquely in this way, \cite[\Prop 2.2, 2.3]{Hav90}, \cite[\Prop 2.1]{Ter92}.

Note that by connectedness assumption, the matrix dimension $D$ (or equivalently any matrix with the same pattern of zero and non-zero entries) is 
indecomposable, \ie, $DD^t$ and $D^tD$ are irreducible square matrices, \cite[\Lem 1.3.2, 2.3.1]{GdHJ89}.

\begin{theorem}\emph{\cite{GiLo19}}.
Let $\N\subset\M$ be as in the previous theorem.

Consider the eigenvalue equations 
\begin{align*}
D^t D \sqrt{\nu} =&\; d^2 \sqrt{\nu}\\
D D^t \sqrt{\mu} =&\; d^2 \sqrt{\mu}
\end{align*}
where $\sqrt{\nu} = (\nu^{1/2}_1,\ldots,\nu^{1/2}_n)^t$, $\sqrt{\mu} = (\mu^{1/2}_1,\ldots,\mu^{1/2}_m)^t$ are vectors with strictly positive entries and $l^2$-normalized (thus unique by Perron-Frobenius theory). 
%$\sum_j \nu_j = 1$ and $\sum_i \mu_i = 1$
Then the minimal expectation $E^0:\M\rightarrow\N$ is determined by
$$(E^0)_{ij} = E^0_{ij}, \quad \lambda^{E^0}_{ij} = \frac{d_{ij}}{d}\frac{\mu_i^{1/2}}{\nu_j^{1/2}} $$
where $E^0_{ij}:\M_{ij} \rightarrow \N_{ij}$ is the unique minimal expectation in each subfactor, if $p_iq_j \neq 0$.
\end{theorem}

As a consequence, we have the following \lqq weighted" additivity formula for the scalar dimension $d$ (thus for $[\M:\N]_0$)
$$d = \sum_{i,j} d_{ij} \nu_j^{1/2} \mu_i^{1/2}.$$
Moreover, by setting $\omega_l (q_j) := \nu_j$ and $\omega_r(p_i) := \mu_i$ we have two faithful states $\omega_l$ and $\omega_r$ on $\Z(\N)$ and $\Z(\M)$ respectively, canonically determined by the inclusion. We call them respectively the \emph{left} and \emph{right state} of $\N\subset\M$.
These states provide a characterization of minimality of $E^0$ which extends the previously mentioned (but not explained) \emph{sphericality condition} in the case of subfactors.
Namely, $E^0$ is the only expectation from $\M$ onto $\N$ fulfilling 
$$\omega_l\circ E^0 = \omega_r\ \circ (E^0)' \quad \text{on} \quad \N'\cap\M$$
where $(E^0)' : \N' \rightarrow \M'$ is the \emph{dual} expectation in the sense of Kosaki \cite{Kos86}. Note that the previous equation makes sense because $E^0(\N'\cap\M) = \Z(\N)$ and $(E^0)'(\N'\cap\M) = \Z(\M)$, and it defines a canonical state on $\N'\cap\M$ for the inclusion, that we call \emph{spherical state} of $\N\subset\M$, denoted by $\omega_s$.

\begin{remark}
In the subfactor case one has that the square root of the minimal index is additive and multiplicative, namely $d = d_1 + d_2$ if $d$ is the dimension of $\N\subset\M$ and $d_1,d_2$ are obtained by cutting with projections $p_1,p_2\in\N'\cap\M$ such that $p_1+p_2 = \oneop$. Moreover, let $\N\subset\M\subset\L$ be two consecutive subfactors, then $d = d_1 d_2$ if $d,d_1,d_2$ are respectively the dimensions of $\N\subset\L$, $\N\subset\M$, $\M\subset\L$. 

These relations do no longer hold for inclusions with non-trivial centers, indeed one has to replace the scalar dimension with the matrix dimension to have $D=D_1+D_2$ and $D=D_2D_1$. In particular, the dimension (thus the minimal index) is in general only submultiplicative $d\leq d_1 d_2$, while the minimal index itself can be additive $d^2 = d_1^2 + d_2^2$ (if one of $\N$ or $\M$ is a factor).
\end{remark}

\section{Extremality and super-extremality}

In the case of finite direct sums of \emph{finite} factors, as for type $\II_1$ subfactors, one can compare the two theories of index (trace/minimal). Given a connected inclusion of such algebras $\N\subset\M$ with finite index, on one hand, we have the matrix dimension $D$, the minimal conditional expectation $E^0$ with index $[\M:\N]_0$, and the spherical state $\omega_s$ on $\N'\cap\M$. On the other hand, we have the Jones index matrix $\Lambda_\N^\M$, the Jones index $[\M:\N]$ and a uniquely determined trace $\tau$ on $\M$, called the \emph{Markov trace} of $\N\subset\M$, \cite[\Sec 2.7,3.7]{GdHJ89}, which extends to the Jones tower.

\begin{definition}
In \cite[\Sec 3]{GiLo19}, we called \emph{extremal} an inclusion which fulfills $E^0 = E^\tau$, where $E^\tau$ is the Markov trace-preserving expectation, and \emph{super-extremal} an inclusion which fulfills in addition $\omega_s = \tau_{\restriction \N'\cap\M}$.
\end{definition}

By \cite[\Prop 3.2]{Hav90}, \cite[\Cor 3.7.4]{GdHJ89} we have $\Ind(E^\tau) = [\M:\N]\oneop$, thus for an extremal inclusion we have $[\M:\N]_0 = [\M:\N]$. Recall from \cite[\Lem 3.2]{GiLo19} that $\omega_s = \tau_{\restriction \N'\cap\M}$ is equivalent to $\omega_l = \tau_{\restriction \Z(\N)}$, where $\omega_l$ is by definition ${\omega_s}_{\restriction{\Z(\N)}}$.
Another condition that one might consider is the equality of matrices $D=\Lambda_{\N}^{\M}$, which corresponds to an \emph{entrywise extremality} for the inclusion once reduced with every $p_iq_j$. Clearly for a finite subfactor, all these notions, including super-extremality, boil down to the ordinary notion of extremality, see Remark \ref{rmk:extremality}.

The following results completely settle the analysis of (super-)extremal inclusions of \emph{multi-matrices} (always assumed to have finite-dimensional centers, thus with finite index). As in the first paragraph of the previous section, if $\N\subset\M$ is the inclusion, denote by $k = (k_1,\ldots,k_n)^t$, $h = (h_1,\ldots,h_m)^t$ the vectors of dimensions such that $\N\cong \bigoplus_{j=1,\ldots,n} M_{k_j}(\CC)$ and $\M\cong \bigoplus_{i=1,\ldots,m} M_{h_i}(\CC)$, and denote by $\Lambda$ the Bratteli inclusion matrix. Recall that $\Lambda k = h$ is the consistency of the Bratteli diagram.

\begin{theorem}\emph{\cite{GiLo19}.}
Let $\N\subset\M$ be a connected multi-matrix inclusion. 

Then $D=\Lambda_\N^\M=\Lambda$ and the inclusion is always extremal, namely $E^0=E^\tau$. The inclusion is also super-extremal if and only if 
$$\Lambda^t h = d^2 k.$$
%\ie, if and only if $h$ (up to normalization) is the Perron-Frobenius eigenvector.. and it implies the same for $k$
\end{theorem}

The index (we need not specify which one) of a super-extremal multi-matrix inclusion is easy to compute and it has the following properties:

\begin{proposition}\emph{\cite{GiLo19}.}
Let $\N\subset\M$ be as in the previous theorem. If the inclusion is super-extremal then
$$[\M:\N]_0 = \frac{\|h\|_{l^2}^2}{\|k\|_{l^2}^2}$$
\emph{(}the ratio of the algebraic dimensions of $\M$ over $\N$\emph{)}. In particular the index must be a positive integer (because rational and algebraic integer) and every positive integer (not only squares of integers) is the index of such an inclusion.

Moreover, the index of super-extremal multi-matrix inclusions is clearly multiplicative.
\end{proposition}

\section{Open problems}

Some natural problems (currently under investigation) that arise from the analysis of the minimal index for von Neumann inclusion reviewed here are:

\begin{problem}
In the case of inclusions with infinite-dimensional and atomic centers there can be more than one expectation such that $\|\Ind(E^0)\| = \inf_{E}\{\|\Ind(E)\|\}$. Can one find a preferred, canonical one, whose index is scalar and which is related to the (infinite) matrix dimension in some way? 
\end{problem}

Note that in the case of the previous problem the matrix dimension can be defined as for finite-dimensional centers using minimal central projections.

\begin{problem}
In the same situation as above, does the theory of minimal index admit a purely 2-$C^*$-categorical (or better 2-$W^*$-categorical) formulation, \cf \cite{GLR85}? Namely, does the theory of intrinsic tensor-categorical dimension admit an extension beyond  tensor $C^*$-categories with simple tensor unit or with finitely reducible tensor unit? What is a \lqq standard solution" of the conjugate equations beyond the previously mentioned regimes?
\end{problem}

We do not want to ask the same question of categorical translation beyond the case of atomic centers, because we cannot immagine by now a good notion of direct integral of objects in a tensor $C^*$- (or $W^*$-) category.

\begin{problem}
In the case of inclusions with infinite-dimensional and possibly diffuse centers (in the absence of minimal central projections) what is a good substitute of the matrix dimension?
%a map between centers
%Can one still find a canonical minimal expectation with scalar index and which is related to $D$ in some way? 
%or do it the other way around trova E^delta e poi calcolaci D con estensioni olomorfe
\end{problem}

\begin{problem}
Study the consequences of and characterize super-extremality for inclusion of finite direct sums of type $\II_1$ factors. Find examples of such inclusions that go beyond tensoring a super-extremal multi-matrix inclusion with a type $\II_1$ factor. 
\end{problem}

\begin{problem}
Study $C^*$-Frobenius algebra objects in tensor \Cstar-categories with non-simple but finitely reducible tensor unit (\eg for unitary multi-fusion categories), \cf \cite[\Ch 3]{BKLR15}, \cite[\Ch 4]{EGNO15}. Study the relation between the realization of $C^*$-Frobenius algebras in $\End(\bigoplus_{i=1}^{n}\R)$ (Q-systems) or $\Bim(\bigoplus_{i=1}^{n}\R)$ and extensions of $\bigoplus_{i=1}^{n}\R$, where $\R$ is a factor.
\end{problem}

It is known that every unitary fusion category, and more generally every rigid tensor \Cstar-category with simple unit can be realized as endomorphisms or bimodules of a factor (that can be chosen either of type $\II_1$ or of type $\III$, and in some cases hyperfinite), \cite{HaYa00}, \cite{Yam03}, \cite{BHP12}, \cite{GiYu19}.

\begin{problem}
Study the problem of realizing unitary multi-fusion categories, or more generally rigid 2-\Cstar-categories with finite-dimensional centers, as endomorphisms or bimodules of $\bigoplus_{i=1}^{n}\R$, where $\R$ is a factor. Is every such abstract category realizable in operator-algebras, if so, is the realization unique in a suitable sense? 
\end{problem}

\bigskip
{\bf Acknowledgements.}
We thank the Institute of Mathematics \lqq Simion Stoilow" of the Romanian Academy and the West University in Timi\c{s}oara for the kind hospitality during the 27th International Conference in Operator Theory (OT27) in July, 2--6, 2018. We thank the organizers for the invitation to present our work at the conference and for financial support. We also thank Maria Stella Adamo and Yoh Tanimoto for discussions.

%\newpage
\small

%%%NOWALLINTEX%%%\bibliography{mybib}%%%Links to mybib.bib same folder

\def\cprime{$'$}

\end{document}